\input amstex
\documentstyle{amsppt}
\document
\topmatter
\title
Complex foliations and K\"ahler QCH surfaces
\endtitle
\author
W\l odzimierz Jelonek
\endauthor

\abstract{The aim of this paper is to describe   complex
2-dimensional foliations in a K\"ahler surface. }
\thanks{MS Classification: 53C55,53C25,53B35. Key words and phrases:
K\"ahler surface,  quasi constant holomorphic sectional
curvature,QCH manifold, homothetic foliation, conformal foliation
}\endthanks
 \endabstract
\endtopmatter
\define\G{\Gamma}

\define\n{\nabla}
\define\om{\omega}
\define\Om{\Omega}

\define\w{\wedge}
\define\k{\diamondsuit}
\define\th{\theta}

\define\a{\alpha}

\define\lb{\lambda}

\define\1{D_{\lb}}
\define\2{D_{\mu}}
\define\0{\Omega}

\define\De{\Cal D}

\define\m{(M,g,J)}
\define \E{\Cal E}
\bigskip
{\bf 1. Introduction. } The aim of the present paper is to
investigate complex 2-dimensional foliations on connected K\"ahler
surfaces $(M,g,J)$.   Two dimensional complex distribution $\De$
on a K\"ahler surface defines an opposite almost Hermitian
structure $I$.  We  prove that if $\De$ is integrable then $I$ is
Hermitian if and only if $\De$ is totally geodesic and
holomorphic.  Also  $I$ is Hermitian if  $\De$ is conformal and
totally geodesic.  Hence any complex totally geodesic and
conformal foliation is holomorphic.  We also prove that the
following conditions are equivalent for a complex,
  2-dimen\-sional  foliation $\De$ on a K\"ahler surface $(M,g,J)$:
(a)  $\De$ is conformal with $L_Vg=\th(V) g$ on $\De^{\perp}$, (b)
$\De$ is quasi holomorphic i.e. $L_{\xi}J(\E)\subset \De$ for
$\xi\in\G(\De)$,
(c)$2g(\n_X\xi,Y)=\th(\xi)g(X,Y)+J\th(\xi)\om(X,Y)$ for
$\xi\in\G(\De)$ and $X,Y\in\E$ where $\om$ is the K\"ahler form of
$(M,g,J)$ and $d\Om_I=2\th\w\Om_I$ where $\Om_I$ is the K\"ahler
form of $(M,g,I)$.

A QCH K\"ahler surface is a K\"ahler surface
 admitting a global, $2$-dimensional, $J$-invariant
distribution $\De$ having the following property: The holomorphic
curvature $K(\pi)=R(X,JX,JX,X)$ of any $J$-invariant $2$-plane
$\pi\subset T_xM$, where $X\in \pi$ and $g(X,X)=1$, depends only
on the point $x$ and the number
$|X_{\De}|=\sqrt{g(X_{\De},X_{\De})}$, where $X_{\De}$ is an
orthogonal projection of $X$ on $\De$.(see [J-1],[J-2],[J-4]
[G-M]). We classify compact, simply connected QCH surfaces  for
which the opposite almost Hermitian structure $I$ induced by $\De$
is Hermitian and one of the distributions  $ker(J\circ I\pm Id)$
(hence $\De$ or $\Cal E=\De^{\perp}$) is integrable.   We also
classify  K\"ahler surfaces $(M,g,J)$ admitting a solution of an
equation $H^\phi=g$ for a certain function $\phi\in C^\infty(M)$
where $H^\phi$ is the Hessian of $\phi$.

\bigskip
{\bf 2.   Complex foliations on a K\"ahler surface.} Our notation
is as in [K-N]. We start with (see  [B-W], [Ch-N]):
\medskip
{\bf Definition 1.} A foliation $\Cal F$ on a Riemannian manifold
$(M,g)$ is called conformal if
$$L_Vg=\a(V)g$$ holds on $T\Cal F^{\perp}$ where $\a$ is a one
form vanishing on $T\Cal F^{\perp}$ for any $V$ tangent to the
leaves od $\Cal F$. A foliation $\Cal F$ is called homothetic if
it is conformal and $d\a=0$.
\medskip
{\bf Definition 2.} A complex distribution $\De$ on a complex
manifold $(M,g,J)$ is called holomorphic if
$L_{\xi}J(TM)\subset\De$  for any $\xi\in\G(\De)$.
\medskip

In the sequel we shall not make the distinction between a
foliation and the distribution of vectors tangent to the leaves of
the foliation.
 A foliation (distribution) $\De$ on
a complex manifold $(M,J)$ is called complex if $J(\De)\subset
\De$.  Complex 2-dimensional foliation (distribution)  on a
complex surface defines an opposite almost Hermitian structure
defined by $I_{|\De}=J_{|\De},I_{|\De^{\perp}}=-J_{|\De^{\perp}}$.
If we denote by $\Om_J,\Om_I$ the K\"ahler forms of $(M,g,J)$ and
$(M,g,I)$ respectively  then  $\Om_J^2=-\Om_I^2$.  In the
following all the considered distributions $\De,\E=\De^{\perp}$
are assumed to be complex and 2-dimensional. We shall recall some
results from [G-M-1]. Let
$$R(X,Y)Z=([\n_X,\n_Y]-\n_{[X,Y]})Z\tag 2.1$$ and let us write $$R(X
,Y,Z,W)=g(R(X,Y)Z,W).$$ If $R$ is the curvature tensor of a QCH
K\"ahler manifold $\m$, then there exist functions $a,b,c\in
C^{\infty}(M)$ such that
$$R=a\Pi+b\Phi+c\Psi,\tag 2.2$$
where $\Pi$ is the standard K\"ahler tensor of constant
holomorphic curvature i.e.
$$\gather \Pi(X,Y,Z,U)=\frac14(g(Y,Z)g(X,U)-g(X,Z)g(Y,U)\tag 2.3\\
+g(JY,Z)g(JX,U)-g(JX,Z)g(JY,U)-2g(JX,Y)g(JZ,U)),\endgather $$ the
tensor $\Phi$ is defined by the following relation
$$\gather \Phi(X,Y,Z,U)=\frac18(g(Y,Z)h(X,U)-g(X,Z)h(Y,U)\tag 2.4\\+g(X,U)h(Y,Z)-g(Y,U)h(X,Z)
+g(JY,Z)h(JX,U)\\-g(JX,Z)h(JY,U)+g(JX,U)h(JY,Z)-g(JY,U)h(JX,Z)\\
-2g(JX,Y)h(JZ,U)-2g(JZ,U)h(JX,Y)),\endgather$$ and finally
$$\Psi(X,Y,Z,U)=-h(JX,Y)h(JZ,U)=-(h_J\otimes h_J)(X,Y,Z,U).\tag 2.5$$
where $h_J(X,Y)=h(JX,Y)$.
\medskip
We say that an almost Hermitian structure $I$ satisfies the second
Gray curvature condition if
$$R(X,Y,Z,W)-R(IX,IY,Z,W)=R(IX,Y,IZ,W)+R(IX,Y,Z,JW).\tag {G2}$$
  Hence
$(M,g,I)$ satisfies the second Gray condition if $I$ preserves the
Ricci tensor and $W^+$ (with respect to the orientation given by
$I$) is degenerate. Note that for a QCH surface $(M,g,J)$ the
distribution $\De$ induces an almost Hermitian structure $I$
satisfying the second Gray condition.
\bigskip
{\bf 3. An almost hermitian structure $I$ induced by complex
foliation $\De$.}  We start with characterizing two dimensional
foliations in 4-dimensional Riemannian manifold.

\bigskip
{\bf Theorem 1.}  {\it  Let $(M,g)$ be a Riemann manifold, dim
$M=4$ and let $\De$ be a two dimensional distribution on $M$ . Let
$\{\th_1,\th_2,\th_3,\th_4\}$ be an orthonormal basis of 1-forms
such that $\De=\{X:\th_3(X)=\th_4(X)=0\}$. Let $\om=\th_3\w\th_4$.
Then $\De$ is integrable if and only if $d\om=\phi\w\om$ for a
certain $\phi\in\Cal A^1(M)$.}
\medskip
{\it Proof}.    $\De$ is integrable if and only if
$d\th_3=a\w\th_3+b\w\th_4,d\th_4=c\w\th_3+d\w\th_4$. It follows
that if $\De$ is integrable then
$$ d\om=d\th_3\w\th_4-\th_3\w d\th_4=a\w\th_3\w\th_4-\th_3\w
d\w\th_4=(a+d)\w\om. $$

Assume now that $d\om=\phi\w\om$.   Let
$d\th_3=a\w\th_3+b\w\th_4+k\th_1\w\th_2,
d\th_4=c\w\th_3+d\w\th_4+l\th_1\w\th_2$. Then
$d\om=(a+d)\w\om+k\th_1\w\th_2\w\th_4-l\th_3\w\th_1\w\th_2$. Hence
$(\phi-a-d)\om=k\th_1\w\th_2\w\th_4-l\th_3\w\th_1\w\th_2$.
Consequently $k=l=0$ and $\De$ is integrable.$\k$

Remark:  Let $(M,g,J)$ be a K\"ahler surface with a K\"ahler form
$\om=\Om_J$. Let $\De$ be a complex foliation  $J\De\subset\De$.
Let us denote by $\Cal E=\De^{\perp}$ the orthogonal complement of
$\De$. Let us define by $I$ an almost Hermitian structure given by
$I_{|\De}=J_{|\De}$, $I_{|\Cal E}=-J_{|\Cal E}$. Then
$\Om_J=\om_1+\om_2,\Om_I=\om_1-\om_2$ where $\om_1\in
\bigwedge^2\De,\om_2\in \bigwedge^2\E$ .We have
$d\om_2=\th\w\om_2$ where $\th^{\sharp}\in\G(\De) $. Thus
$d\Om_I=-2\th\w\Om_I$.

\bigskip

{\bf Theorem 2.} {\it Let  $(M,g,J)$ be a K\"ahler surface and
$\De$ be a complex foliation, dim$\De=2$. Then the almost
Hermitian structure $I$ given by
$I_{|\De}=J_{|\De},I_{|\De^{\perp}}=-J_{|\De^{\perp}}$ is
integrable if and only if   $\De$ is holomorphic and totally
geodesic.}

\medskip
{\it Proof.}  Let us assume that $I$ is integrable.   Let
$\om_2=\th_3\w\th_4$ where $\De=ker\th_3\cap\th_4$.  Then
$d\om_2=\th\w\om_2$ where $\th\in \bigwedge^1\De$.   If
$\Om_J=\om_1+\om_2$ then $\Om_I=\om_1-\om_2$ and
$d\Om_I=-2\th\w\om_2=2\th\w\Om_I$. Since $\th^{\sharp}\in\De$ and
$I$ is integrable we have $\n_XI=0$ for
$X\in\De=span\{\th,J\th\}$. It follows that $\De$ is totally
geodesic.    In fact $IX=JX$ for $X,Y\in\De$ implies
$I(\n_YX)=J(\n_YX)$.  We show that $\De$ is holomorphic. We have
to show that $g(\n_{JY}X-J\n_YX,Z)=0$ for arbitrary $Y$ and
$Z\in\Cal E$.  If $Y$ belongs to $\De$ it is a consequence of the
fact that  $\De$ is totally geodesic. Let $Y\in \E$. Then
$g(\n_{JX}\xi,Y)-g(J\n_X\xi,Y)=g(\xi,\n_{JX}Y+\n_XJY)$ for
 $\xi\in\De,X,Y\in\Cal E$. We have to show that for $X,Y\in \Cal
E$ we obtain  $\n_{JX}Y+\n_XJY\in\Cal E$.  Taking $Z=JY$ it is
equivalent to  $-\n_{JX}JZ+\n_XZ\in\Cal E$ for $X,Z\in\Cal E$.
Since $IZ=-JZ$ to $\n_XIZ+I(\n_XZ)=-J(\n_XZ)$.  Hence
$\n_XIZ=-2J(\n_XZ_{|\De})$. But $I$ is integrable, hence
$\n_{IX}IIZ=-2J(\n_{IX}IZ_{|\De})=-2J(\n_XZ_{|\De})$. It implies
$-\n_{JX}JZ+\n_XZ\in\Cal E$.  Assume now that $\De$ is totally
geodesic and holomorphic.  Then $\n_XI=0$ for $X\in\De$.
 If $X,Y\in\De,Z\in\E$ then
$\n_XY\in\De,\n_XZ\in\E$. It implies $\n_XIY=\n_XIZ=0$. Since
$\De$ is holomorphic we get $g(\n_{JX}\xi-J\n_X\xi,Y)=0$ for
$Y\in\E$ and arbitrary $X$. Let now $Y\in\De,X\in\E$. Then
$\n_XIY=(J-I)(\n_XY)=2J(\n_XY_{|\E})$. Note that
$\n_{IX}IIY=-2J(\n_{JX}JY_{|\E})=2J(\n_XY_{|\E})$. Now assume that
$X,Y\in\E$. Then $\n_XIY=(-I-J)(\n_XY)=-2J(\n_XY_{|\De})$. Again
using that $\De$ is holomorphic we obtain
$g(\xi,\n_{JX}Y)+g(\xi,\n_XJY)=0$ and
consequently$g(\xi,\n_XY-\n_{JX}JY)=0$.  Thus
$\n_XY_{|\De}=\n_{JX}JY_{|\De}$ i $\n_{IX}IIY=\n_XIY$ also in this
case .$\k$

\bigskip
{\bf Theorem 3.} {\it Let $(M,g,J)$ be a K\"ahler surface. Let us
assume that $\De$ is a complex foliation , $dim\De=2$ and the
structure $I$ defined by $\De$ is integrable.  Then

$$2g(\n_X\xi,Y)=\th(\xi)g(X,Y)+J\th(\xi)\om(X,Y)\tag *$$
for every $\xi\in\G(\De)$ i $X,Y\in \Cal E$.  On the other hand if
$\De$ is totally geodesic and  (*) is satisfied  then $I$ is
integrable.}
\medskip
{\it Proof.}  We have $d\om_2=\th\w\om_2$ and
$d\Om_I=-2\th\w\Om_I$. Note that $\n_XI=0$ for $X\in\De$. Indeed
in the interior  of the set  $F=\{x:\th_x=0\}$ the structure $I$
is K\"ahler and outside  $F$ we have $span
\{\th^\sharp,J\th^\sharp\}=\De$. Since $I$ is Hermitian we have
$\n_{\th}I=\n_{J\th}I=0$.  Let $\xi\in\G(\De)$. Then
$L_{\xi}-\n_{\xi}=-\n\xi$.  Note that $\n_\xi\om_1=\n_\xi\om_2=0$.
It implies $(L_{\xi}-\n_\xi)\om_2=L_\xi\om_2=\xi\lrcorner
d\om_2=\xi\lrcorner(\th\w\om_2)=\th(\xi)\om_2$.  Hence $
g(J\n_X\xi,Y)+g(JX,\n_Y\xi)=\th(\xi)\om_2(X,Y)$ and
 $
g(\n_X\zeta,Y)-g(X,\n_Y\zeta)=-\th(J\zeta)g(JX,Y)$.  Since $I$ is
Hermitian $\De$ is holomorphic and
$g(\n_{JX}\xi,Y)+g(JX,\n_Y\xi)=\th(\xi)g(JX,Y)$. Thus

$g(\n_X\xi,Y)+g(X,\n_Y\xi)=\th(\xi)g(X,Y)$ for any $X,Y\in\E$.
Consequently
$$2g(\n_X\xi,Y)=\th(\xi)g(X,Y)+J\th(\xi)\om(X,Y).$$   On the other
hand if $\De$ is totally geodesic and  (*) is satisfied then $\De$
is holomorphic and  $I$ is integrable.$\k$

\medskip
{\bf Theorem 4.} {\it Let $\De$ be a 2-dimensional, complex
foliation on a
 K\"ahler surface $(M,g,J)$.  Then
$\E=\De^{\perp}$ is integrable if and only if $I$ is almost
K\"ahler.}

\medskip
{\it Proof.}  If $\E,\De$ are integrable then
$d\om_1=\phi\w\om_1,d\om_2=\psi\w\om_2$ and $d\om_1=-d\om_2$ hence
$d\om_1=d\om_2=0$. On the other hand if $d\om_1=d\om_2=0$ then
$\De,\E$ are integrable.$\k$

\bigskip
{\bf Theorem 5.} {\it Let  $(M,g,J)$ be a K\"ahler surface and
$\De$ be a complex distribution $J\De\subset\De$ and $I$ be an
almost Hermitian  structure defined by $\De$. Then $\De\subset
ker\n I=\{X:\n_XI=0\}$ if and only if   $\De$ is totally
geodesic.}
\medskip
{\it Proof}  $\Rightarrow$  Let  $X\in\De$ then  $IX=JX$ and for
$Y\in\De$ we get  $I\n_YX=J\n_YX$  and consequently $\n_YX\in\De$.

 $\Leftarrow$  Let  $\De$ be totally geodesic. Then
 $\n_XY\in\De(\E)$ if $Y\in\De(\E)$ and $X\in\De$.  It implies
 $\n_XIY+I\n_XY=\pm J\n_XY$ and $\n_XIY=0$ for
 $Y\in\De(Y\in\E)$.$\k$
\medskip
{\it Remark.}  Let $(M,g,J)$ be a 3-symmetric K\"ahler surface
(see [K]). Then   $(M,g,J)$ is a QCH surface and it admits
strictly almost K\"ahler opposite structure  $I$.   If we denote
by $\De$ the K\"ahler nullity of $I$ then $(M,g,J)$ is a QCH with
respect to $\De$ and $\De$ is totally geodesic.  Hence there exist
totally geodesic complex foliations on K\"ahler surfaces for which
the almost Hermitian structure $I$ induced by $\De$ is not
integrable.

\medskip
{\bf Lemma.  } {\it Let us assume that  $J\De=\De$ and
$\E=\De^{\perp}$. If
$\om_1\in\bigwedge^2\De,\om_2\in\bigwedge^2\E$ and
$\om=\om_1+\om_2$ where $\om$ is a  K\"ahler form of $(M,g,J)$
then $\n_X\om_2(Y,Z)=0$ for $X\in\De,Y,Z\in\E$.}

\medskip
{\it Proof  }  We have  $\om_2(Y,Z)=g(JY,Z)$. Note also that
$\om_2(Y,Z)=g(JY,Z)$ if $Z\in\E$ and $Y$ is arbitrary. Thus
$$\gather \n_X\om_2(Y,Z)+\om_2(\n_XY,Z)+\om_2(Y,\n_XZ)=g(J\n_XY,Z)+g(JY,\n_XZ)\\
=\om_2(\n_XY,Z)+\om_2(Y,\n_XZ)\endgather$$$\k$
\medskip

{\bf Corollary 1.}  {\it  If  $\De$  is complex and integrable
then
$$  g(\n_X\xi,Y)-g(\n_Y\xi,X)=J\th(\xi)g(JX,Y)$$
for $\xi\in\G(\De)$ and $X,Y\in\E$  where  $d\Om_I=-2\th\w\Om_I$.
For $X,Y\in\G(\E)$ we have $[X,Y]_{|\De}=-J\th^\sharp\om(X,Y)$. }

\medskip
{\it Proof.}  We have  $\Om_J=\om_1+\om_2,\Om_I=\om_1-\om_2$ and
$d\om_2=\th\w\om_2$ where $\th^\sharp\in\G(\De)$.  Hence
$(L_\xi-\n_\xi)\om_2=
L_\xi\om_2-\n_\xi\om_2=\th(\xi)\om_2-\n_\xi\om_2$. It follows that
for $Y,Z\in\E$ we get
$\om_2(\n_X\xi,Y)+\om_2(X,\n_Y\xi)=\th(\xi)\om_2(X,Y)$ and
$g(\n_XJ\xi, Y)-g(X,\n_YJ\xi)=\th(\xi)g(JX,Y)$. Thus $g(\n_X\zeta,
Y)-g(X,\n_Y\zeta)=J\th(\zeta)g(JX,Y)$.

Hence   $(**)[X,Y]_{|\De}=-J\th^\sharp\om(X,Y)$.$\k$
\medskip
{\it Remark}  Let us note that  $d\th(X,Y)=0$ for $X,Y\in\De$ and
$d\th(X,Y)=0$ for $X,Y\in\E$.  In fact from (**) it follows for
$X,Y\in\E$. Next $d\th\w(\th_1\w\th_2-\th_3\w\th_4)=0$. Hence the
result follows.

\bigskip
 {\bf Theorem  6.}  {\it The following conditions are equivalent for a complex,
  2-dimen\-sional  foliation $\De$ on a K\"ahler surface $(M,g,J)$:

(a)  $\De$ is conformal

(b) $\De$ is quasi holomorphic  i.e.  $L_{\xi}J(\E)\subset \De$
for $\xi\in\G(\De)$

(c)$2g(\n_X\xi,Y)=\th(\xi)g(X,Y)+J\th(\xi)\om(X,Y)$ for
$\xi\in\G(\De)$ and $X,Y\in\E$ where $d\Om_I=-2\th\w\Om_I$.}
\medskip
{\it Proof.}   $(a)\Rightarrow(b),(c)$  Note that
$g(\n_X\xi,Y)-g(\n_Y\xi,X)=J\th(\xi)g(JX,Y)$ and
$g(\n_X\xi,Y)+g(\n_Y\xi,X)=\phi(\xi)g(X,Y)$.  Summing up
$$2g(\n_X\xi,Y)=J\th(\xi)g(JX,Y)+\phi(\xi)g(X,Y).$$  Thus
  $g(\n_{JX}\xi,JY)=g(\n_X\xi,Y)$.  It is also easily seen that
  $\phi=\th$.

$(b)\Rightarrow(a)$  If   $\De$ is quasi holomorphic then from
$$g(\n_X\xi,JY)-g(\n_{JY}\xi,X)= J\th(\xi)g(JX,JY)$$ we get
$$-g(\n_XJ\xi,Y)-g(\n_YJ\xi,X)=-\th(J\xi)g(X,Y)$$ and consequently
$$g(\n_X\xi,Y)+g(\n_Y\xi,X)=\th(\xi)g(X,Y)$$ and $\De$ is
conformal.

\medskip
{\bf Corollary 2.}  {\it   If  $\De$ is a complex foliation on a
K\"ahler surface  $(M,g,J)$ which is conformal and totally
geodesic then the opposite almost hermitian structure is
Hermitian.  Hence  conformal and totally geodesic foliation $\De$
is holomorphic. }

\medskip
{\it Proof.}  Since $\De$ is conformal is quasi holomorphic. Since
it is totally geodesic it is holomorphic.  Hence the result
follows from Th.2.$\k$

\medskip
{\bf  Theorem  7.}  {\it  Let  $(M,g,J)$ be a compact simply
connected QCH K\"ahler surface. If $I$ is Hermitian and one of the
distributions $ker(I\circ J\pm Id)$ is integrable and  then
 $(M,g,J)$ is
globally of Calabi type or a product of Riemann surfaces. }

{\it Proof.}  Let  $\De=ker(J\circ I-Id)$ be integrable. Then
$\De$ is totally geodesic homothetic foliation. In fact since $I$
is integrable we have $d\th^-=0$ and since $M$ is compact $d\th=0$
(see [A-G]).  Note that  that homothetic foliations are classified
locally  in [Ch-N].  Thus if we do not assume that $M$ is simply
connected but is compact it follows that locally is isometric to a
Calabi type surface or a product of Riemann surfaces. In the
compact simply connected case the result follows from [J-3].$\k$

{\bf Theorem 8.} {\it  Let $\De$ be a complex totally geodesic and
holomorphic foliation on a K\"ahler surface $(M,g,J)$. Then $\De$
is conformal, $2g(\n_X\xi,Y)=|\th|^2g(X,Y)$ where $\th=g(\xi,.)$
and $d\overline{\om}=-2\th\w\overline{\om}$,
$L_{\eta}g=\th(\eta)g$ on $\De^{\perp}$ for $\eta\in\G(\De)$. If
$d\th^-=0$ then $(M,g,J)$ is a  QCH K\"ahler surface  whose
negative almost Hermitian  $G_2$ structure $\overline{J}$ is
Hermitian. On the other hand if $(M,g,J)$ is a QCH K\"ahler
surface with integrable distribution  $\De$ and Hermitian
structure $ \overline{J}$ then $\De$ is a totally geodesic
conformal holomorphic foliation with $d\th^-=0$.}$\k$

\medskip

{\it Proof}  We have for $X,Y\in\G(\De^{\perp})$
$2g(\n_X\xi,Y)=|\th|^2g(X,Y)$ thus for $Z\in\G(\De^{\perp})$ we
obtain
$$\gather 2g(\n_Z\n_X\xi,Y)+2g(\n_X\xi,\n_ZY)=\\Z|\th|^2g(X,Y)+|\th|^2g(\n_ZX,Y)+|\th|^2g(X,\n_ZY)\endgather$$
If $d\th^-=0$ then $d\th(JX,JY)=-d\th(X,Y)$.   Hence
$$d\th(\xi,X)=-\frac12X|\th|^2, d\th(J\xi,X)=-\frac12JX|\th|^2$$
and
$\n_ZY=-\frac12g(Z,Y)\xi-\frac12\om(Z,Y)J\xi+\n_ZY_{\De^{\perp}}$.
Consequently
$$\gather
2g(\n_Z\n_X\xi,Y)=\frac12g(Z,Y)X|\th|^2+\frac12\om(Z,Y)JX|\th|^2\\+Z|\th|^2g(X,Y)+|\th|^2g(\n_ZX,Y)\endgather$$
Hence
$$\gather 4g(R(Z,X)\xi,Y)=-X|\th|^2g(Z,Y)+Z|\th|^2g(X,Y)+\\\om(Z,Y)JX|\th|^2-\om(X,Y)JZ|\th|^2=0.\endgather$$

It follows that   $(M,g,\overline{J})$ is a $G_3$  Gray manifold,
and since $\overline{J}$ is Hermitian it follows that
$(M,g,\overline{J})$ is a $G_2$ manifold and $(M,g,J)$ is a  QCH
manifold. We have assumed here that $\th\ne0$ but if $\th=0$ then
 $(M,g,J)$ is a product hence also a  $QCH$ manifold.$\k$

Now we consider  K\"ahler QCH surfaces which are semi-symmetric
(i.e. $R.R=0$) with Hermitian structure $\overline{J}$ and such
that $(M,g,\overline{J})$ is of Gray type $(G_1)$. These are
semi-symmetric surfaces foliated by two dimensional Euclidean
space (see [J-5]). Let $\om\in\bigwedge^2\De$ be as above. We have
$\n_{\xi}\om=0$ and $-2\n_X\om=-JX\w\th+X\w J\th$ for $X\in
\De^{\perp}$.  Hence $\n_X\n_\xi\om=0$,
$$-2\n_\xi\n_X\om=-J(\n_\xi X)\w\th-JX\w\n_\xi\th+\n_\xi X\w
J\th+X\w J\n_\xi\th. $$

Note that
$\n_X\xi=\frac12|\th|^2X+\frac12X\ln|\th|^2\xi+\frac12JX\ln|\th|^2J\xi$,
and
$$\gather -2\n_{[X,\xi]}\om=-J[X,\xi]_\perp\w\th+[X,\xi]_\perp\w
J\th=\\-J\n_X\xi_\perp\w\th+J\n_\xi X\w\th+\frac12|\th|^2X\w
J\th-\n_\xi X\w J\th=\\-\frac12|\th|^2 JX\w\th+J\n_\xi X\w\th
+\frac12|\th|^2X\w J\th-\n_\xi X\w J\th\endgather$$

Hence    $$ -2R(\xi,X).\om=JX\w\n_\xi\th-X\w
J\n_\xi\th+\frac12|\th|^2 JX\w\th-\frac12|\th|^2X\w J\th$$

Consequently (note that $R(\xi,X)=0$)
$\n_\xi\th+\frac12|\th|^2\xi=0$ i $\n_\xi\xi=-\frac12|\th|^2\xi$
and $\n_{J\xi}\xi=\frac12|\th|^2J\xi$. We also have
$\xi|\th|^2=-|\th|^4$ and $[\xi,J\xi]=-|\th|^2J\xi$

Now we shall show that $\eta=\frac1{|\th|^2}\xi$ is a holomorphic
vector field and $J\eta$ is a holomorphic divergence free vector
field.

Note that
$\n_X\eta=\frac12X-\frac12X\ln|\th|^2\eta+\frac12JX\ln|\th|^2J\eta$
for $X\in\De^{\perp}$.  Hence $\n_{JX}\eta=J\n_X\eta$.

Similarly   $\n_{J\xi}\eta=\frac12J\xi$ and
$\n_\xi\eta=\frac12\xi$. Hence $\n_X\eta=\frac12 X$ for $X\in\De$.
It is also clear that $\n_{JX}\eta=J\n_{X}\eta$ if $X\in\De$ and
$\eta\in\frak{hol}(M)$. We also have
$\n_XJ\eta=\frac12JX-\frac12X\ln|\th|^2J\eta-\frac12JX\ln|\th|^2\eta$
for $X\in\De^{\perp}$ and $\n_XJ\eta=\frac12JX$ for $X\in\De$.
Consequently $div J\eta=0$.

If $d\th=0$ then   $\th=-d\ln|\th|^2 $ and consequently on
semi-symmetric surfaces foliated by two dimensional Euclidean
space of Calabi type we have $\th\ne0$ on the whole of $M$.

{\bf  Theorem 9.}  {\it   Let $(M,g,J)$ be a K\"ahler surface.
Then the following conditions are equivalent:

(1) There exists a vector field $\xi$ such that $\n\xi=c I,
c\in\Bbb R-\{0\}$

(2)There exists a vector field $\eta$ such that $\n\eta=c J,
c\in\Bbb R-\{0\}$

(3) There exists a function $\phi\in C^{\infty}(M)$ such that
$H^\phi=c g$ where $H^\phi$ is the Hessian of $\phi$ and $c\in
\Bbb R-\{0\}$,

(4) $(U,g,J)$ is a semi-symmetric K\"ahler surface of Calabi type
where $U$ is an open dense subset of $M$ ,

(5)There exists an open and dense subset $U$ of $M$ such that
$(U,g,J)$ is locally isometric to

$$ g=zg_{\Sigma}+\frac1{Cz}dz^2+Cz(dt+\alpha)^2$$
where $(\Sigma,g_{\Sigma})$ is a Riemannian surface with area form
$\om_{\Sigma}$ and $d\alpha=\om_{\Sigma}$. The K\"ahler form of
$(U,g,J)$ is $\Om=z\om_{\Sigma}+dz\w(dt+\alpha)$. }

\medskip
{\it Proof.} Let us assume that (1) holds. It is easy to see that
(1),(2)  are equivalent and (3) implies (1).

Note that  $\n_XJ\xi=cJX$ so we can take  $\eta=J\xi$. Let us
define the distribution $\De=span\{\xi,J\xi\}$. We can assume that
$c=\frac12$. It is clear that $\xi$ is a holomorphic vector field
and $J\xi$ is a holomorphic
 Killing vector field.  Note that $\xi,J\xi$ are different from zero on an open dense subset $U$ of $M$.
   What is more if $T=\n J\xi=\frac12J$
 then
$R(X,\xi)Y=\n T(X,Y)=0$. Hence $R(X,Y)\xi=R(X,Y)J\xi=0$ and
$$R(X,JX,JX,X)=||X_{\De^\perp}||^4K(\De^{\perp})$$ where
$K(\De^{\perp})$ is a sectional curvature of the distribution
$\De^{\perp}$. It follows that $(M,g,J)$ is a QCH K\"ahler surface
i $R=c\Psi$ and $R.R=0$. The distribution $\De$ is totally
geodesic in particular is integrable. Since
$L_{J\xi}g=0,L_{\xi}g=g$ it follows that $\De$ is a complex
conformal foliation and the almost Hermitian structure $I$
determined by  $\De$ is Hermitian. Note that
$L_{\eta}g=\th(\eta)g$ on $\De^{\perp}$ for $\eta\in\G(\De)$ where
$\th^{\sharp}=\frac1{|\xi|^2}\xi$.  Hence $|\th|=\frac1{|\xi|}$.
Vector field $\xi=\frac1{|\th|^2}\th^{\sharp}$ is holomorphic. One
can easily verify that $d\th=0$ since $X|\th|^2=0$ for
$X\in\De^\perp$ since
$\n_X\xi=\frac12X-\frac12X\ln|\th|^2\xi+\frac12JX\ln|\th|^2J\xi$
for $X\in\De^{\perp}$.  Hence   $(M,g,J)$ is a  semi-symmetric
K\"ahler surface of Calabi  type.  Note that since
$\th=-d\ln|\th|^2$  on $U$ it follows that if the distribution
$\De$ extends over the whole of $M$ and consequently $\th$ is
defined on $M$ then $\th\ne0$ on the whole of $M$ and consequently
$U=M$. Note that the function $\phi=\frac1{|\th|^2}$ satisfies (3)
and the field $\xi=\n\phi$ satisfies (1).

On the other hand if  $(M,g,J)$ is a  semi-symmetric  K\"ahler
surface of Calabi type then vector field  $\xi
=\frac1{|\th|^2}\th^\sharp$ satisfies $\n_X\xi=\frac12X$.
Semi-symmetric K\"ahler surfaces of Calabi type are classified in
[J-5] which gives equivalence $(4)\Leftrightarrow (5)$.$\k$

{\it Remark.}  It is known ( see [T] )that the only complete
K\"ahler surface satisfying (3)  is a Euclidean space $(\Bbb
C^2,can)$ with standard metric $can$. Let $\{e_1,e_2,e_3,e_4\}$ be
a standard orthonormal basis of $\Bbb C^2$, $Je_1=e_2,Je_3=e_4$.
Then $\xi= x_1e_1+y_1e_2+x_2e_3+y_2e_4,
J\xi=x_1e_2-y_1e_1+x_2e_4-y_2e_3$ where $z_1=x_1+iy_1,z_2+iy_2$
are standard complex coordinates on $\Bbb C^2$ and
$\phi=\frac12(x_1^2+y_1^2+x_2^2+y_2^2)$.  In this case we have
$U=\Bbb C^2-\{0\}$ and the totally geodesic complex foliation
$\De=span\{\xi,J\xi\}$  defines on $U$ a Hermitian conformally
K\"ahler  structure $I$ which does not extend to the whole of
$\Bbb C^2$ and $\xi(0)=0$. Hence $(\Bbb C^2-\{0\},can)$ is a
semi-symmetric surface of Calabi type  (clearly $R=0$ in this
case). Note that orthogonal complex structures on domains of $\Bbb
C^2$ are described in [S-V].

\bigskip
\centerline{\bf References.}

\par
\medskip
\cite{A-G} V. Apostolov, P. Gauduchon { The Riemannian
Goldberg-Sachs Theorem}  Internat. J. Math. vol.8, No.4,
(1997),421-439

\par
\medskip
\cite{B-W} P. Baird and J. Wood {\it Harmonic morphisms between
Riemannian manifolds} London Mathematical Society Monographs 29
(Oxford University Press, Oxford 2003)
\par
\medskip
\cite{Bes}  A. L. Besse {\it Einstein manifolds}, Ergebnisse,
ser.3, vol. 10, Springer-Verlag, Berlin-Heidelberg-New York, 1987.
\par
\medskip
\cite{Ch-N} S.G.Chiossi and P-A. Nagy {\it Complex homothetic
foliations on K\"ahler manifolds}  Bull. London Math. Soc. 44
(2012) 113-124.

\par
\medskip
\cite{G-M} G.Ganchev, V. Mihova {\it K\"ahler manifolds of
quasi-constant holomorphic sectional curvatures}, Cent. Eur. J.
Math. 6(1),(2008), 43-75.

\par
\medskip
\cite{J-1} W. Jelonek,{ Compact holomorphically pseudosymmetric
K\"ahler manifolds} Coll. Math.117,(2009),No.2,243-249.
\par
\medskip
\cite {J-2} W.Jelonek {\it K\"ahler manifolds with quasi-constant
holomorphic curvature}, Ann. Glob. Anal. and Geom, vol.36, p.
143-159,( 2009)
\par
\medskip
\cite{J-3} W. Jelonek, {\it K\"ahler manifolds with homothetic
foliation by curves} arxiv
\par
\medskip
\cite{J-4} W.  Jelonek {\it K\"ahler surfaces with quasi-constant
holomorphic curvature} to appear in Glasgow Math.J.
\par
\medskip
\cite{J-5} W.  Jelonek {\it Semi-symmetric K\"ahler surfaces}
arxiv

\par
\medskip
\cite{K-N} S. Kobayashi and K. Nomizu {\it Foundations of
Differential Geometry}, vol.2, Interscience, New York  1963
\par
\medskip
\cite{K} O. Kowalski {\it Generalized symmetric spaces} Lecture
Notes in Math. 805, Springer, New York,1980.
\par
\medskip
\cite{S-V} Salamon S.,Viaclovsky J., {\it Orthogonal complex
structures on domains in $\Bbb R^4$} Math. Ann.343, 853-899 (2009)
\par
\medskip
\cite{T} Y. Tashiro, Complete Riemannian manifolds and some vector
fields, Trans. Amer. Math. Soc. 117(1965) 251– 275
\par
\medskip

\par
\medskip
Institute of Mathematics

Cracow University of Technology

Warszawska 24

31-155      Krak\'ow,  POLAND.

 E-mail address: wjelon\@pk.edu.pl
\bigskip

\enddocument